\documentclass[graybox]{svmult}

\usepackage{mathptmx}       % selects Times Roman as basic font
\usepackage{helvet}         % selects Helvetica as sans-serif font
\usepackage{courier}        % selects Courier as typewriter font
\usepackage{type1cm}        % activate if the above 3 fonts are
\usepackage{makeidx}         % allows index generation
\usepackage{graphicx}        % standard LaTeX graphics tool
\usepackage{multicol}        % used for the two-column index
\usepackage[bottom]{footmisc}% places footnotes at page bottom

\usepackage{url}

\usepackage{amsmath}
\usepackage{amssymb}
\makeindex             % used for the subject index
\newcommand{\B} {{\mathbb B}}
\newcommand{\Hmat} {{\mathbb H}}
\newcommand{\R} {\mathbb R}
\newcommand{\id} {\mathbb I}
\newcommand{\U} {\mathbb U}
\newcommand{\V} {\mathbb V}
\newcommand{\W} {\mathbb W}
\newcommand{\veps}{\varepsilon}
\newcommand{\vt}{\vartheta}
\newcommand{\cancel}[1]{}
\newcommand{\raw}{\rightarrow}

\newcommand{\fourvec}[4]{\begin{bmatrix} #1 \\ #2 \\ #3 \\ #4 \end{bmatrix}}
\newcommand{\fourvecr}[4]{\left[\begin{array}{r}  #1 \\ #2 \\ #3 \\ #4 \end{array}\right] }

\newcommand{\Hmt}{{[\psi]}}
\newcommand{\Hmttwo}{{[\psi^2]}}
\newcommand{\Hmtthree}{{[\psi^3]}}
\newcommand{\Bern}{{[\beta]}}
\newcommand{\BernTwo}{{[\beta^2]}}
\newcommand{\BernThree}{{[\beta^3]}}
\newcommand{\BernStwo}{{[\xi^2]}}
\newcommand{\BernStwoprime}{{[\xi^2]^{[0,1]}}}
\newcommand{\BernSthree}{{[\xi^3]}}
\newcommand{\BernSthreeprime}{{[\xi^3]^{[0,1]}}}
\newcommand{\HermStwo}{{[\vartheta^2]}}
\newcommand{\HermStwoprime}{{[\vartheta^2]^{[0,1]}}}
\newcommand{\HermSthree}{{[\vartheta^3]}}
\newcommand{\HermSthreeprime}{{[\vartheta^3]^{[0,1]}}}
\newcommand{\sldeg}{\textnormal{sldeg}}

\DeclareMathAlphabet{\mathcal}{OMS}{cmsy}{m}{n}
\newcommand{\cP}{{\mathcal P}}
\newcommand{\cQ}{{\mathcal Q}}
\newcommand{\cS}{{\mathcal S}}

\begin{document}

\title*{Hermite and Bernstein Style Basis Functions for Cubic Serendipity Spaces on Squares and Cubes}
\titlerunning{Hermite and Bernstein Style Basis Functions for Cubic Serendipity Spaces}
\author{Andrew Gillette}
\institute{Department of Mathematics, University of Arizona, 617 N. Santa Rita Ave., P.O. Box 210089, Tucson, AZ 85721, \email{agillette@math.arizona.edu}
}
\maketitle

\abstract{We introduce new Hermite-style and Bernstein-style geometric decompositions of the cubic serendipity finite element spaces $\cS_3(I^2)$ and $\cS_3(I^3)$,  as defined in the recent work of Arnold and Awanou [\emph{Found. Comput. Math.} \textbf{11} (2011), 337--344].
The serendipity spaces are substantially smaller in dimension than the more commonly used bicubic and tricubic Hermite tensor product spaces - 12 instead of 16 for the square and 32 instead of 64 for the cube - yet are still guaranteed to obtain cubic order \textit{a priori} error estimates in $H^1$ norm when used in finite element methods.
The basis functions we define have a canonical relationship both to the finite element degrees of freedom as well as to the geometry of their graphs; this means the bases may be suitable for applications employing \textit{isogeometric analysis} where domain geometry and functions supported on the domain are described by the same basis functions.
Moreover, the basis functions are linear combinations of the commonly used bicubic and tricubic polynomial Bernstein or Hermite basis functions, allowing their rapid incorporation into existing finite element codes.
}

\section{Introduction}

Serendipity spaces offer a rigorous means to reduce the degrees of freedom associated to a finite element method while still ensuring optimal order convergence.
The `serendipity' moniker came from the observation of this phenomenon among finite element practitioners before its mathematical justification was fully understood; see e.g.~\cite{Ci02,H1987,M1990,SF73}.
Recent work by Arnold and Awanou~\cite{AA2011,AA2012} classifies serendipity spaces on cubical meshes in $n\geq 2$ dimensions by giving a simple and precise definition of a space of polynomials $\cS_r(I^n)$ that must be spanned, as well as a unisolvent set of degrees of freedom for them.
Crucially, the space $\cS_r(I^n)$ contains all polynomials in $n$ variables of total degree at most $r$, a property shared by the space of polynomials $\cQ_r(I^n)$ spanned by the standard order $r$ tensor product method.
This property allows the derivation of an \textit{a priori} error estimate for serendipity methods of the same order (with respect to the width of a mesh element) as their standard tensor product counterparts.

In this paper, we provide two coordinate-independent geometric decompositions for both $\cS_3(I^2)$ and $\cS_3(I^3)$, the cubic serendipity spaces in two and three dimensions, respectively.
More precisely, we present sets of polynomial basis functions, prove that they provide a basis for the corresponding cubic serendipity space, and relate them canonically to the domain geometry.
Each basis is designated as either Bernstein or Hermite style, as each function restricts to one of these common basis function types on each edge of the square or cube.
The standard pictures for $\cS_3(I^2)$ and $\cS_3(I^3)$ serendipity elements, shown on the right of Figures~\ref{fg:sq-bicubic} and~\ref{fg:cube-tricubic} below, have one dot for each vertex and two dots for each edge of the square or cube.
We refer to these as \textbf{domain points} and will present a canonical relationship between the defined bases and the domain points.

\begin{figure}[ht]
\begin{tabular}{ccc}
\parbox{.32\textwidth}{\includegraphics[width=.32\textwidth]{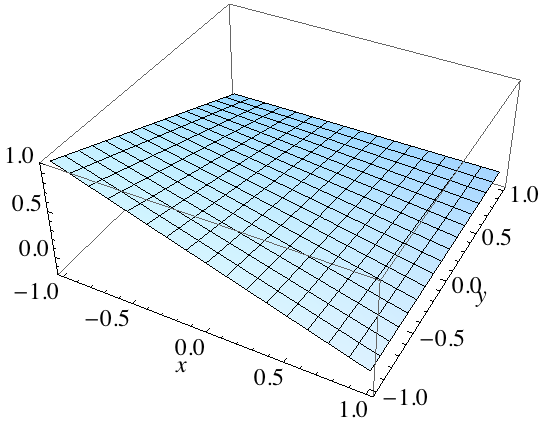}} &
\parbox{.32\textwidth}{\includegraphics[width=.32\textwidth]{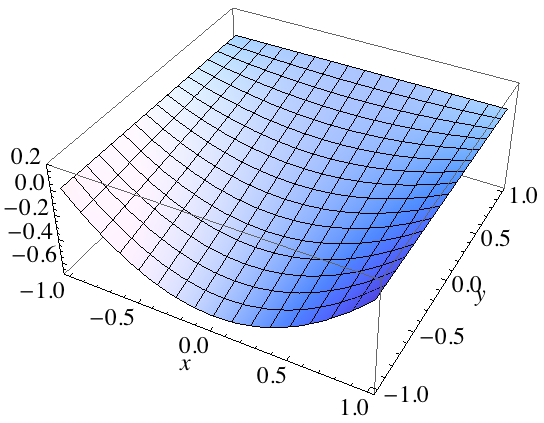}} &
\parbox{.32\textwidth}{\includegraphics[width=.32\textwidth]{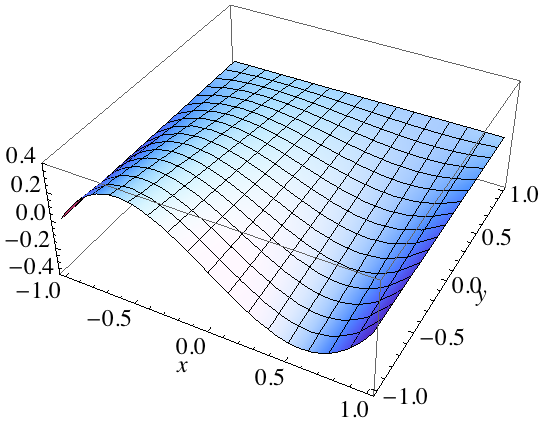}} \\
{$\frac{1}{4} (x-1) (y-1)$} & {$-\frac{1}{4} \sqrt{\frac{3}{2}} \left(x^2-1\right) (y-1)$} & {$-\frac{1}{4} \sqrt{\frac{5}{2}} x \left(x^2-1\right) (y-1)$}
\end{tabular}
\caption{Cubic serendipity functions on $I^2$ from \cite{SB1991}. 
The left function is associated to the vertex below the peak.
The middle and right functions are associated to the edge $y=-1$ but do not correspond to the domain points $(\pm\frac 13,-1)$ in any canonical or symmetric fashion, making them less useful for geometric modeling or isogeometric analysis.}
\label{fg:sb-functions}
\end{figure}

To the author's knowledge, the only basis functions previously available for cubic serendipity finite element purposes employ Legendre polynomials, which lack a clear relationship to the domain points.
Definitions of these basis functions can be found in Szab{\'o} and Babu{\v{s}}ka~\cite[Sections 6.1 and 13.3]{SB1991}; the two functions from~\cite{SB1991} associated to the edge $y=-1$ of $I^2$, are shown in Figure~\ref{fg:sb-functions} (middle and right).
The restriction of these functions to the edge gives an even polynomial in one case and an odd polynomial in the other, forcing an \textit{ad hoc} choice of how to associate the functions to the corresponding domain points $(\pm\frac 13,-1)$.
The functions presented in this paper do have a natural correspondence to the domain points of the geometry.

Maintaining a concrete and canonical relationship between domain points and basis functions is an essential component of the growing field of \textit{isogeometric analysis} (IGA).
One of the main goals of IGA is to employ basis functions that can be used both for geometry modeling and finite element analysis, exactly as we provide here for cubic serendipity spaces.
Each function is a linear combination of bicubic or tricubic Bernstein or Hermite polynomials; the specific coefficients of the combination are given in the proofs of the theorems.
This makes the incorporation of the functions into a variety of existing application contexts relatively easy.
Note that tensor product bases in two and three dimensions are commonly available in finite element software packages (e.g.\ deal.II~\cite{BHK2007}) and cubic tensor products in particular are commonly used both in modern theory (e.g.\ isogeometric analysis~\cite{EH2011}) and applications (e.g.\ cardiac electrophysiology models~\cite{Zetal2012}).
Hence, a variety of areas of computational science could directly employ the new cubic serendipity basis functions presented here.

The benefit of serendipity finite element methods is a significant reduction in the computational effort required for optimal order (in this case, cubic) convergence.
Cubic serendipity methods on meshes of squares requires 12 functions per element, an improvement over the 16 functions per element required for  bicubic tensor product methods.
On meshes of cubes, the cubic serendipity method requires 32 functions per element instead of the 64 functions per element required for tricubic tensor product methods.
Using fewer basis functions per element reduces the size of the overall linear system that must be solved, thereby saving computational time and effort.
An additional computational advantage occurs when the functions presented here are used in an isogeometric fashion.
The process of converting between computational geometry bases and finite element bases is a well-known computational bottleneck in engineering applications~\cite{CHB2009} but is easily avoided when basis functions suited to both purposes are employed.

The outline of the paper is as follows.
In Section~\ref{sec:bkgd}, we fix notation and summarize relevant background on Bernstein and Hermite basis functions as well as serendipity spaces.
In Section~\ref{sec:2d-bases}, we present polynomial Bernstein and Hermite style basis functions for $\cS_3(I^2)$ that agree with the standard bicubics on edges of $I^2$ and provide a novel geometric decomposition of the space.
In Section~\ref{sec:3d-bases}, we present polynomial Bernstein and Hermite style basis functions for $\cS_3(I^3)$ that agree with the standard tricubics on edges of $I^3$, reduce to our bases for $I^2$ on faces of $I^3$, and provide a novel geometric decomposition of the space.
Finally, we state our conclusions and discuss future directions in Section~\ref{sec:conc}.

\section{Background and Notation}
\label{sec:bkgd}

\subsection{Serendipity Elements}
\label{subsec:ser-bkgd}
We first review the definition of serendipity spaces and their accompanying notation from the work of Arnold and Awanou~\cite{AA2011,AA2012}.

\begin{definition}
\label{def:superlin}
The \textbf{superlinear degree} of a monomial in $n$ variables, denoted $\sldeg(\cdot)$, is given by
\begin{equation}
\label{eq:superlin}
\sldeg(x_1^{e_1}x_2^{e_2}\cdots x_n^{e_n}) := \left(\sum_{i=1}^n e_i\right) - \#\left\{e_i~:~e_i=1\right\}
\end{equation}
\end{definition}
In words, $\sldeg(q)$ is the ordinary degree of $q$, ignoring variables that enter linearly.
For instance, the superlinear degree of $xy^2z^3$ is 5.
\begin{definition}
\label{def:srdpty}
Define the following spaces of polynomials, each of which is restricted to the domain $I^n=[-1,1]^n\subset\R^n$:
\begin{align*}
\cP_r(I^n) & := \text{span}_{\R}\left\{\text{monomials in $n$ variables with total degree at most $r$}\right\} \\
\cS_r(I^n) & := \text{span}_{\R}\left\{\text{monomials in $n$ variables with superlinear degree at most $r$}\right\} \\
\cQ_r(I^n) & := \text{span}_{\R}\left\{\text{monomials in $n$ variables of degree at most $r$ in each variable}\right\} %\\
\end{align*}
\end{definition}
Note that $\cP_r(I^n)\subset \cS_r(I^n)\subset \cQ_r(I^n)$, with proper containments when $r$, $n>1$.
The space $\cS_r(I^n)$ is called the degree $r$ \textbf{serendipity space} on the $n$-dimensional cube $I^n$.
In the notation of the more recent paper by Arnold and Awanou~\cite{AA2012}, the serendipity spaces discussed in this work would be denoted $\cS_r\Lambda^0(I^n)$, indicating that they are differential 0-form spaces.
The space $\cQ_r(I^n)$ is associated with standard tensor product finite element methods; the fact that $\cS_r(I^n)$ satisfies the containments above is one of the key features allowing it to retain an $O(h^r)$ \textit{a priori} error estimate in $H^1$ norm, where $h$ denotes the width of a mesh element~\cite{BS2002}.
The spaces have dimension given by the following formulas~(cf.~\cite{AA2011}).
\begin{align*}
\dim \cP_r(I^n) & = {n+r \choose n} \\
\dim \cS_r(I^n) & = \sum_{d=0}^{\min(n,\lfloor r/2\rfloor)}2^{n-d}{n\choose d}{r-d\choose d}\\
\dim \cQ_r(I^n) & = (r+1)^n
\end{align*}
We write out standard bases for these spaces more precisely in the cubic cases of concern here.
\begin{align}
\cP_3(I^2) & = \text{span}\{~1~~,~~\underbrace{x,y}_{\text{linear}}~~,~~\underbrace{x^2,y^2,xy}_{\text{quadratic}}~~,~~\underbrace{x^3,y^3,x^2y,xy^2}_{\text{cubic}}~\} \label{def:p3-2d}\\
\cS_3(I^2) & = \cP_3(I^2)\cup \text{span}\{\underbrace{x^3y,xy^3}_{\text{superlinear cubic }}\} \label{def:s3-2d}\\
\cQ_3(I^2) & = \cS_3(I^2)\cup \text{span}\{x^2y^2,x^3y^2,x^2y^3,x^3y^3\} \label{def:q3-2d} 
\end{align}
Observe that the dimensions of the three spaces are 10, 12, and 16, respectively.
\begin{align}
\cP_3(I^3) & = \text{span}\{1,\underbrace{x,y,z}_{\text{linear}}, \underbrace{x^2,y^2,z^2,xy,xz,yz}_{\text{quadratic}},\underbrace{x^3,y^3,z^3,x^2y,x^2z,xy^2,y^2z,xz^2,yz^2,xyz}_{\text{cubic}}\} \label{def:p3-3d}\\
\cS_3(I^3) & = \cP_3(I^3)\cup \text{span} \{\underbrace{x^3y,x^3z,y^3z,xy^3,xz^3,yz^3,x^2yz,xy^2z,xyz^2, x^3yz,xy^3z,xyz^3}_{\text{superlinear cubic}}\} \label{def:s3-3d}\\
\cQ_3(I^3) & = \cS_3(I^3)\cup \text{span} \{x^3y^2,\ldots,x^3y^3z^3\} \label{def:q3-3d}
\end{align}
Observe that the dimensions of the three spaces are 20, 32, and 64, respectively.

The serendipity spaces are associated to specific \textbf{degrees of freedom} in the classical finite element sense.
For a face $f$ of $I^n$ of dimension $d\geq 0$, the degrees of freedom associated to $f$ for $\cS_r(I^n)$ are (cf.~\cite{AA2011})
\[u\longmapsto\int_f uq,\qquad q\in\cP_{r-2d}(f).\]
For the cases considered in this work, $n=2$ or 3 and $r=3$, so the only non-zero degrees of freedom are when $f$ is a vertex ($d=0$) or an edge ($d=1$).
Thus, the degrees of freedom for our cases are the values
\begin{equation}
\label{eq:srdp-dofs}
u(v),\quad \int_e u~dt,\quad\text{and}\quad\int_e ut~dt,
\end{equation}
for each vertex $v$ and each edge $e$ of the square or cube.

\subsection{Cubic Bernstein and Hermite Bases}

For cubic order approximation on square or cubical grids, tensor product bases are typically built from one of two alternative bases for $\cP_3([0,1])$:
\[\Bern =  \fourvec{\beta_1}{\beta_2}{\beta_3}{\beta_4} := 
%  \fourvecr{1-3x+3x^2-x^3}{x-2x^2+x^3}{x^2-x^3}{x^3} =  
   \fourvec{(1-x)^3}{(1-x)^2x}{(1-x)x^2}{x^3}
\qquad\qquad
\Hmt  =  \fourvec{\psi_1}{\psi_2}{\psi_3}{\psi_4}  :=
  \fourvecr {1 - 3x^2+2x^3}{x-2x^2+x^3}{x^2-x^3}{3x^2-2x^3}
  \]
The set $\{\beta_1,3\beta_2,3\beta_3,\beta_4\}$ is the \textbf{cubic Bernstein} basis and the set $\Hmt$ is the \textbf{cubic Hermite} basis.
Bernstein functions have been used recently to provide a geometric decomposition of finite element spaces over simplices~\cite{AFW2009}. 
Hermite functions, while more common in geometric modeling contexts~\cite{M2006} have also been studied in finite element contexts for some time~\cite{CR1972}.
The Hermite functions have the following important property relating them to the geometry of the graph of their associated interpolant:
\begin{equation}
\label{eq:hmt-interp-01}
u = u(0)\psi_1 + u'(0)\psi_2 - u'(1)\psi_3 + u(1)\psi_4,\qquad\forall u\in\cP_3([0,1]).
\end{equation}
We have chosen these sign and basis ordering conventions so that both bases have the same symmetry property:
\begin{equation}
\label{eq:b-and-h-sym}
\beta_k(1-x)=\beta_{5-k}(x),\qquad \psi_k(1-x)=\psi_{5-k}(x).
\end{equation}
The bases $\Bern$ and $\Hmt$ are related by $\Bern=\V\Hmt$ and $\Hmt=\V^{-1}\Bern$ where 
\begin{equation}
\label{eq:def-V}
\V= \begin{bmatrix} 
1 &-3 & 0 & 0 \\
0 & 1 & 0 & 0 \\
0 & 0 & 1 & 0 \\
0 & 0 &-3 & 1 
\end{bmatrix}
,\qquad
\V^{-1}= \begin{bmatrix} 
1 & 3 & 0 & 0 \\
0 & 1 & 0 & 0 \\
0 & 0 & 1 & 0 \\
0 & 0 & 3 & 1 
\end{bmatrix}.
\end{equation}

Let $[\beta^n]$ denote the tensor product of $n$ copies of $\Bern$.
Denote $\beta_i(x)\beta_j(y)\in[\beta^2]$ by $\beta_{ij}$ and $\beta_i(x)\beta_j(y)\beta_k(z)\in[\beta^3]$ by $\beta_{ijk}$.
In general, $[\beta^n]$ is a basis for $\cQ_3([0,1]^n)$, but we will make use of the specific linear combination used to prove this, as stated in the following proposition.

\begin{proposition}
\label{prop:hmtB-prec}
For $0\leq r,s,t\leq 3$, the reproduction properties of $\Bern$, $\BernTwo$, and $\BernThree$ take on the respective forms
\begin{align}
x^r &= \sum_{i=1}^4 {3-r\choose 4-i}\beta_i,\label{eq:hmtB-prec} \\
x^ry^s &= \sum_{i=1}^4\sum_{j=1}^4 {3-r\choose 4-i}{3-s\choose 4-j}\beta_{ij},
\label{eq:bi-hmtB-prec} \\
x^ry^sz^t &= \sum_{i=1}^4\sum_{j=1}^4\sum_{k=1}^4 {3-r\choose 4-i}{3-s\choose 4-j}{3-t\choose 4-k}\beta_{ijk}.\label{eq:tri-hmtB-prec}
\end{align}
\end{proposition}

The proof is elementary.
We have a similar property for tensor products of the Hermite basis $\Hmt$, using analogous notation.
The proof is a simple matter of swapping the order of summation.

\begin{proposition}
\label{prop:hmt-prec}
Let
\begin{equation}
\label{eq:def-veps}
\veps_{r,i}:=\sum_{a=1}^4{3-r\choose 4-a}v_{ai}
\end{equation}
where $v_{ai}$ denotes the $(a,i)$ entry (row, column) of $\V$ from (\ref{eq:def-V}).
For $0\leq r,s,t\leq 3$, the reproduction properties of $\Hmt$, $\Hmttwo$, and $\Hmtthree$ take on the respective forms
\begin{align}
x^r &= \sum_{i=1}^4 \veps_{r,i}\psi_i,\label{eq:hmt-prec} \\
x^ry^s &= \sum_{i=1}^4\sum_{j=1}^4 \veps_{r,i}\veps_{s,j}\psi_{ij},
\label{eq:bi-hmt-prec} \\
x^ry^sz^t &= \sum_{i=1}^4\sum_{j=1}^4\sum_{k=1}^4 \veps_{r,i}\veps_{s,j}\veps_{t,k}\psi_{ijk}.
\label{eq:tri-hmt-prec}
\end{align}
\end{proposition}

Transforming the bases $\Bern$ and $\Hmt$ to domains other than $[0,1]$ is straightforward.
If $T:[a,b]\raw[0,1]$ is linear, then replacing $x$ with $T(x)$ in each basis function expression for $\Bern$ and $\Hmt$ gives bases for $\cP_3([a,b])$.
Note, however, that the derivative interpolation property for $\Hmt$ must be adjusted to account for the scaling:
\begin{align}
u(x) =  u(a)\psi_1(T(x)) & + (b-a) u'(a)\psi_2(T(x)) \notag \\
 & - (b-a)  u'(b)\psi_3(T(x)) + u(b)\psi_4(T(x)),\quad\forall u\in\cP_3([a,b]).
 \label{eq:hmt-interp-I}
\end{align}
In geometric modeling applications, the coefficient $(b-a)$ is sometimes left as an adjustable parameter, usually denoted $s$ for scale factor~
\cite{HLS1993},
however, $(b-a)$ is the only choice of scale factor that allows the representation of $u$ given in (\ref{eq:hmt-interp-I}).
For all the Hermite and Hermite style functions, we will use \textbf{derivative-preserving scaling} which will include scale factors on those functions related to derivatives; this will be made explicit in the various contexts where it is relevant.

\begin{remark}
\label{rmk:dofs}
Both $\Bern$ and $\Hmt$ are Lagrange-like at the endpoints of $[0,1]$, i.e.\ at an endpoint, the only basis function with non-zero value is the function associated to that endpoint ($\beta_1$ or $\psi_1$ for 0, $\beta_4$ or $\psi_4$ for 1).
This means the two remaining basis functions of each type ($\beta_2$, $\beta_3$ or $\psi_2$, $\psi_3$) are naturally associated to the two edge degrees of freedom (\ref{eq:srdp-dofs}).
We will refer to these associations between basis functions and geometrical objects as the standard \textbf{geometrical decompositions} of $\Bern$ and $\Hmt$.
\end{remark}

\section{Local Bases for $\mathcal{S}_3(I^2)$}
\label{sec:2d-bases}

\begin{figure}[h]
\centering
\includegraphics[scale=.35]{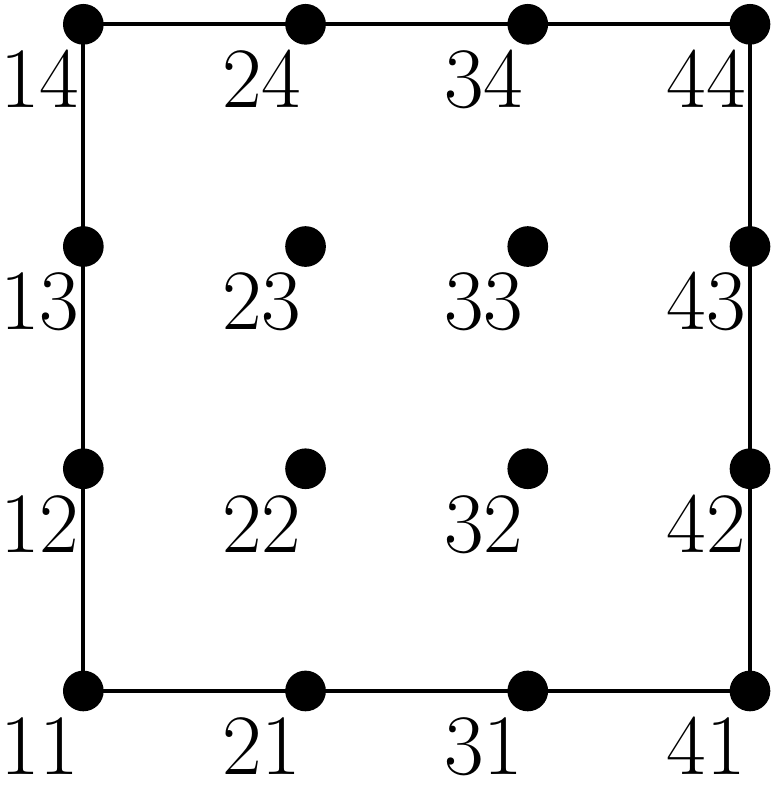}\qquad\qquad
\includegraphics[scale=.35]{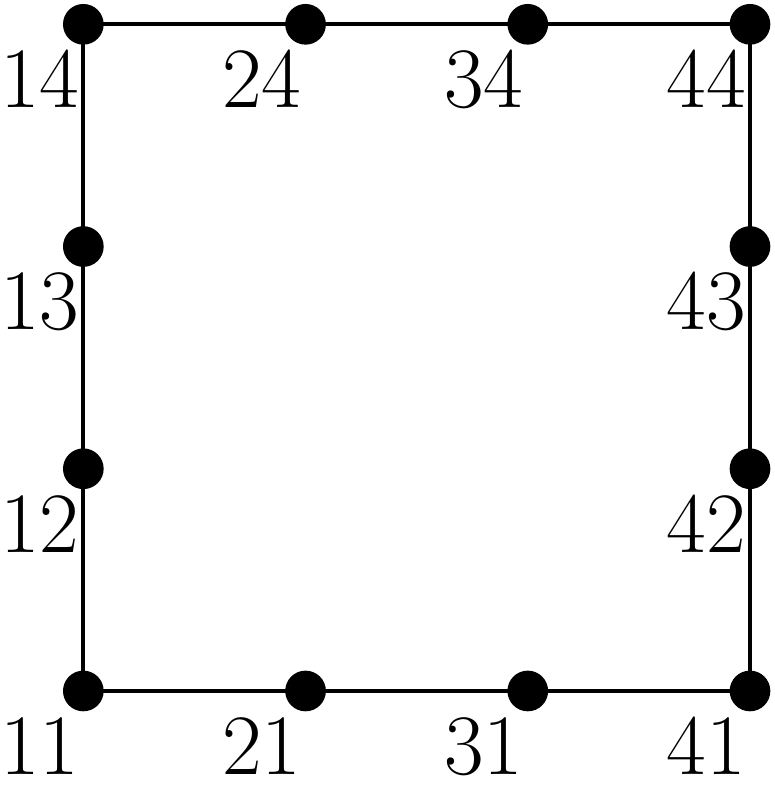}
\caption{On the left, ordered pairs from $X$ are shown next to the domain point of $[0,1]^2$ to which they correspond.  On the right, only those ordered pairs used for the serendipity basis are shown.  The correspondences $V\leftrightarrow$ vertices, $E\leftrightarrow$ edge points, and $D\leftrightarrow$ domain interior points are evident. }
\label{fg:sq-bicubic}
\end{figure}

Before defining local bases on the square, we fix notation for the domain points to which they are associated.
For $[0,1]^2$, define the set of ordered pairs
\[
X := \left\{ \{i,j\} ~~|~~ i,j\in \{1, \ldots, 4\} \right\}.
\]
Then $X$ is the disjoint union $V\cup E\cup D$ where
\begin{align}
V & := \left\{ \{i,j\}\in X ~~|~~ i,j\in \{1, 4\} \right\}; \label{eq:sq-dom-V}\\
E & := \left\{ \{i,j\}\in X ~~|~~ \text{exactly one of $i,j$ is an element of }\{1, 4\}\right\}; \label{eq:sq-dom-E}\\
D & := \left\{ \{i,j\}\in X ~~|~~ i,j\in \{2, 3\} \right\} \label{eq:sq-dom-D}
\end{align}
The $V$ indices are associated with vertices of $[0,1]^2$, the $E$ indices to edges of $[0,1]^2$, and the $D$ vertices to the domain interior to $[0,1]^2$.
The relation between indices and domain points of the square is shown in Figure~\ref{fg:sq-bicubic}.
We will frequently denote an index set $\{i,j\}$ as $ij$ to reduce notational clutter.

\subsection{A local Bernstein style basis for $\mathcal{S}_3(I^2)$}
\label{subsec:2d-bern}

We now establish a local Bernstein style basis for $\mathcal{S}_3(I^3)$ where $I:=[-1,1]$.
Define the following set of 12 functions, indexed by ${V\cup E}$; note the scaling by $1/16$.
\begin{align}
\BernStwo =
\begin{bmatrix}
\xi_{11} \\
\xi_{14} \\
\xi_{41}  \\
\xi_{44} \\
\xi_{12} \\
\xi_{13} \\
\xi_{42} \\
\xi_{43} \\
\xi_{21} \\
\xi_{31} \\
\xi_{24} \\
\xi_{34}
\end{bmatrix}
 & =
\begin{bmatrix}
(1-x) (1-y) (-2-2 x+x^2-2 y+y^2) \\
(1-x) (y+1) (-2-2 x+x^2+2 y+y^2) \\
(x+1) (1-y) (-2+2 x+x^2-2 y+y^2) \\
(x+1) (y+1) (-2+2 x+x^2+2 y+y^2) \\
(1-x) (1-y)^2 (y+1) \\
(1-x) (1-y) (y+1)^2 \\
(x+1) (1-y)^2 (y+1) \\
(x+1) (1-y) (y+1)^2 \\
(1-x)^2 (x+1) (1-y) \\
(1-x) (x+1)^2 (1-y) \\
(1-x)^2 (x+1) (y+1) \\
(1-x) (x+1)^2 (y+1) 
\end{bmatrix}
\cdot \frac 1{16}
\label{eq:bern-ser-2d-def}
%
% Here's the basis on [0,1]^2, simplified
%
%\begin{bmatrix}
% (x-1)(y-1)  (-2 x+x^2+(y-1)^2) \\
% -(x-1) y    (-2 x+x^2+y^2) \\
% -x (y-1)    (x^2-2y+y^2) \\
%   x y       (-1+x^2+y^2) \\
% -(x-1) (y-1)^2 y \\
% (x-1) (y-1) y^2 \\
% x (y-1)^2 y \\
% -x (y-1) y^2 \\
% -(x-1)^2 x (y-1) \\
% (x-1) x^2 (y-1) \\
% (x-1)^2 x y \\
% -(x-1) x^2 y
%\end{bmatrix}
%
% Here's the basis on [0,1]^2, expanded
%
%  =
%\begin{bmatrix}
% 1-3 x+3 x^2-x^3-3 y+5 x y-3 x^2 y+x^3 y+3 y^2-3 x y^2-y^3+x y^3 \\
% -2 x y+3 x^2 y-x^3 y+y^3-x y^3 \\
% x^3-2 x y-x^3 y+3 x y^2-x y^3 \\
% -x y+x^3 y+x y^3 \\
% y-x y-2 y^2+2 x y^2+y^3-x y^3 \\
% y^2-x y^2-y^3+x y^3 \\
% x y-2 x y^2+x y^3 \\
% x y^2-x y^3 \\
% x-2 x^2+x^3-x y+2 x^2 y-x^3 y \\
% x^2-x^3-x^2 y+x^3 y \\
% x y-2 x^2 y+x^3 y \\
% x^2 y-x^3 y
% \end{bmatrix}
\end{align}
Fix the basis orderings
\begin{align}
\BernStwo &:= [~\underbrace{\xi_{11},\xi_{14},\xi_{41} ,\xi_{44}}_{\text{indices in $V$}},~~\underbrace{\xi_{12},\xi_{13},\xi_{42},\xi_{43},\xi_{21},\xi_{31},\xi_{24},\xi_{34}}_{\text{indices in $E$}}~],\label{eq:order-xi-two}\\
\BernTwo & := [~\underbrace{\beta_{11},\beta_{14},\beta_{41} ,\beta_{44}}_{\text{indices in $V$}},~~\underbrace{\beta_{12},\beta_{13},\beta_{42},\beta_{43},\beta_{21},\beta_{31},\beta_{24},\beta_{34}}_{\text{indices in $E$}},~~\underbrace{\beta_{22},\beta_{23},\beta_{32} ,\beta_{33}}_{\text{indices in $D$}}~] \label{eq:order-beta-two}
\end{align}

%%  Bernstein figure
%\begin{figure}[h]
%\begin{tabular}{ccc}
%\parbox{.32\textwidth}{\includegraphics[width=.32\textwidth]{figs/beta11.png}} &
%\parbox{.32\textwidth}{\includegraphics[width=.32\textwidth]{figs/beta21.png}} &
%\parbox{.32\textwidth}{\includegraphics[width=.32\textwidth]{figs/beta31.png}} 
%\\
%{$\beta_{11}^I$} & {$\beta_{21}^I$} & {$\beta_{31}^I$}
%\\
%\parbox{.32\textwidth}{\includegraphics[width=.32\textwidth]{figs/xi11.png}} &
%\parbox{.32\textwidth}{\includegraphics[width=.32\textwidth]{figs/xi21.png}} &
%\parbox{.32\textwidth}{\includegraphics[width=.32\textwidth]{figs/xi31.png}} 
%\\
%{$\xi_{11}$} & {$\xi_{21}$} & {$\xi_{31}$}
%\end{tabular}
%\caption{The top row shows 3 of the 16 bicubic Bernstein functions on $I^2$ while the bottom row shows 3 of the 12 cubic Bernstein style serendipity functions.  The visual differences are subtle, although some changes in concavity can be observed.  Note that functions in the same column have the same values on the edges of $I^2$.}
%\label{fg:bern-fns}
%\end{figure}

The following theorem will show that $\BernStwo$ is a geometric decomposition of $\cS_3(I^2)$, by which we mean that each function in $\BernStwo$ has a natural association to a specific degree of freedom, i.e.\ to a specific domain point of the element.  

\begin{theorem}
\label{thm:bern-2d}
Let $\beta_{\ell m}^I$ denote the scaling of $\beta_{\ell m}$ to $I^2$, i.e.\
\[\beta_{\ell m}^I :=\beta_{\ell}((x+1)/2)\beta_m((y+1)/2).\]
The set $\BernStwo$ has the following properties:
\begin{enumerate}
\renewcommand{\labelenumi}{(\roman{enumi})}
\item $\BernStwo$ is a basis for $\cS_3(I^2)$.
\item For any $\ell m\in V\cup E$,\quad $\xi_{\ell m}$ is identical to $\beta_{\ell m}^I$ on the edges of $I^2$.
\item $\BernStwo$ is a geometric decomposition of $\cS_3(I^2)$.
\end{enumerate}
\end{theorem}

\begin{proof}
\smartqed
For (i), we scale $\BernStwo$ to $[0,1]^2$ to take advantage of a simple characterization of the reproduction properties.
Let $\BernStwoprime$ denote the set of scaled basis functions $\xi_{\ell m}^{[0,1]}(x,y):=\xi_{\ell m}(2x-1,2y-1)$.
Given the basis orderings in (\ref{eq:order-xi-two})-(\ref{eq:order-beta-two}), it can be confirmed directly that $\BernStwoprime$ is related to $\BernTwo$ by
\begin{equation}
\label{eq:bi-bernS-AbernB}
\BernStwoprime = \B\BernTwo
\end{equation}
where $\B$ is the $12\times 16$ matrix with the structure
\begin{equation}
\label{eq:B-struc}
\B:= \left[ 
\begin{array}{c|c} ~\id~ & ~\B'~ \end{array}
\right],
\end{equation}
where $\id$ is the $12\times 12$ identity matrix and $\B'$ is the $12\times 4$ matrix
\begin{equation}
\label{eq:def-B-prime}
\B' =\begin{bmatrix}
 -4 & -2 & -2 & -1 \\
 -2 & -4 & -1 & -2 \\
 -2 & -1 & -4 & -2 \\
 -1 & -2 & -2 & -4 \\
 2 & 0 & 1 & 0 \\
 0 & 2 & 0 & 1 \\
 1 & 0 & 2 & 0 \\
 0 & 1 & 0 & 2 \\
 2 & 1 & 0 & 0 \\
 0 & 0 & 2 & 1 \\
 1 & 2 & 0 & 0 \\
 0 & 0 & 1 & 2 \\
\end{bmatrix}
.
\end{equation}
Using $ij\in X$ to denote an index for $\beta_{ij}$ and $\ell m\in V\cup E$ to denote an index for $\xi_{\ell m}^{[0,1]}$, the entries of $\B$ can be denoted by $b^{\ell m}_{ij}$ so that
\begin{equation}
\label{eq:B-entries}
\B:= \left[ 
\begin{array}{ccccc}
b^{11}_{11} & \cdots & b^{11}_{ij} & \cdots &  b^{11}_{33} \\ 
\vdots & \ddots & \vdots & \ddots  & \vdots \\
b^{\ell m}_{11} & \cdots & b^{\ell m}_{ij} & \cdots &  b^{\ell m}_{33} \\ 
\vdots & \ddots & \vdots & \ddots  & \vdots \\
b^{34}_{11} & \cdots & b^{34}_{ij} & \cdots &  b^{34}_{33}
\end{array}
\right].
\end{equation}
We now observe that for each $ij\in X$,
\begin{equation}
\label{eq:A-coeff-constr}
{3-r\choose 4-i}{3-s\choose 4-j} = \sum_{\ell m\in V\cup E} {3-r\choose 4-\ell}{3-s\choose 4-m} b^{\ell m}_{ij},
\end{equation}
for all $(r,s)$ pairs such that $\sldeg(x^ry^s)\leq 3$ (recall Definition~\ref{def:superlin}).
Note that this claim holds trivially for the first 12 columns of $\B$, i.e.\ for those $ij\in V\cup E\subset X$.
For $ij\in D\subset X$, (\ref{eq:A-coeff-constr}) defines an invertible linear system of 12 equations with 12 unknowns whose solution is the $ij$ column of $\B'$; the 12 $(r,s)$ pairs correspond to the exponents of $x$ and $y$ in the basis ordering of $\cS_3(I^2)$ given in (\ref{def:p3-2d})-(\ref{def:s3-2d}).
Substituting (\ref{eq:A-coeff-constr}) into (\ref{eq:bi-hmtB-prec}) yields:
\begin{align*}
x^ry^s & = \sum_{ij\in X} \left(\sum_{\ell m\in V\cup E} {3-r\choose 4-\ell}{3-s\choose 4-m} b^{\ell m}_{ij}\right) \beta_{ij}.
\end{align*}
Swapping the order of summation and regrouping yields
\begin{align*}
x^ry^s & = \sum_{\ell m\in V\cup E} {3-r\choose 4-\ell}{3-s\choose 4-m} \left(\sum_{ij\in X} b^{\ell m}_{ij}\beta_{ij} \right).
\end{align*}
The inner summation is exactly $\xi_{\ell m}^{[0,1]}$ by (\ref{eq:bi-bernS-AbernB}), implying that
\begin{equation}
\label{eq:bi-bernS-prec}
x^ry^s = \sum_{\ell m\in V\cup E} {3-r\choose 4-\ell}{3-s\choose 4-m}\xi_{\ell m}^{[0,1]},
\end{equation}
for all $(r,s)$ pairs with $\sldeg(x^ry^s)\leq 3$.
Since $\BernStwoprime$ has 12 elements which span the 12 dimensional space $\cS_3([0,1]^2)$, it is  a basis for $\cS_3([0,1]^2)$.
By scaling, $\BernStwo$ is a basis for $\cS_3(I^2)$.

For (ii), note that an edge of $[0,1]^2$ is described by an equation of the form $\{x~\text{or}~y\}=\{0~\text{or}~1\}$.
Since $\beta_2(t)$ and $\beta_3(t)$ are equal to 0 at $t=0$ and $t=1$, $\beta_{ij}\equiv 0$ on the edges of $[0,1]^2$ for any $ij\in D$.
By the structure of $\B$ from (\ref{eq:B-struc}), we see that for any $\ell m\in V\cup E$,
\begin{equation}
\label{eq:xi-alt-expr}
\xi_{\ell m}^{[0,1]}=\beta_{\ell m}+\sum_{ij\in D}b^{\ell m}_{ij}\beta_{ij}.
\end{equation}
Thus, on the edges of $[0,1]^2$, $\xi_{\ell m}^{[0,1]}$ and $\beta_{\ell m}$ are identical.
After scaling back, we have $\xi_{\ell m}$ and $\beta_{\ell m}^I$ identical on the edges of $I^2$, as desired.

For (iii), the geometric decomposition is given by the indices of the basis functions, i.e.\ the function $\xi_{\ell m}$ is associated to the domain point for $\ell m\in V\cup E$.
This follows immediately from (ii), the fact that $\BernTwo$ is a tensor product basis, and Remark~\ref{rmk:dofs}.
\qed
\end{proof}

\begin{remark}
\label{rmk:proof-tech}
It is worth noting that the basis $\BernStwo$ was derived by essentially the reverse order of the proof of part (i) of the theorem.
More precisely, the twelve coefficients in each column of $\B$ define an invertible linear system given by (\ref{eq:A-coeff-constr}).
After solving for the coefficients, we can immediately derive the basis functions via (\ref{eq:bi-bernS-AbernB}).
By the nature of this approach, the edge agreement property (ii) is guaranteed by the symmetry properties of the basis $\Bern$.
This technique was inspired by a previous work for Lagrange-like quadratic serendipity elements on convex polygons~\cite{RGB2011a}.
\end{remark}

\subsection{A local Hermite style basis for $\mathcal{S}_3(I^2)$}
\label{subsec:2d-hmt}

We now establish a local Hermite style basis $\HermStwo$ for $\cS_3(I^2)$ using the bicubic Hermite basis $\Hmttwo$ for $\cQ_3([0,1]^2)$.
Define the following set of 12 functions, indexed by ${V\cup E}$; note the scaling by $1/8$.
\begin{align}
\HermStwo =
\begin{bmatrix}
\vt_{11} \\
\vt_{14} \\
\vt_{41}  \\
\vt_{44} \\
\vt_{12} \\
\vt_{13} \\
\vt_{42} \\
\vt_{43} \\
\vt_{21} \\
\vt_{31} \\
\vt_{24} \\
\vt_{34}
\end{bmatrix}
 & =
\begin{bmatrix}
-(1-x) (1-y) (-2+x+x^2+y+y^2) \\
-(1-x) (y+1) (-2+x+x^2-y+y^2) \\
-(x+1) (1-y) (-2-x+x^2+y+y^2) \\
-(x+1) (y+1) (-2-x+x^2-y+y^2) \\
 (1-x) (1-y)^2 (y+1) \\
 (1-x) (1-y) (y+1)^2 \\
 (x+1) (1-y)^2 (y+1) \\
 (x+1) (1-y) (y+1)^2 \\
 (1-x)^2 (x+1) (1-y) \\
 (1-x) (x+1)^2 (1-y) \\
 (1-x)^2 (x+1) (y+1) \\
 (1-x) (x+1)^2 (y+1) 
\end{bmatrix}
\cdot \frac 1{8}
\label{eq:hmt-ser-2d-def}
\end{align}

Fix the basis orderings
\begin{align}
\HermStwo &:= [~\underbrace{\vt_{11},\vt_{14},\vt_{41} ,\vt_{44}}_{\text{indices in $V$}},~~\underbrace{\vt_{12},\vt_{13},\vt_{42},\vt_{43},\vt_{21},\vt_{31},\vt_{24},\vt_{34}}_{\text{indices in $E$}}~],\label{eq:order-vt-two}\\
\Hmttwo & := [~\underbrace{\psi_{11},\psi_{14},\psi_{41} ,\psi_{44}}_{\text{indices in $V$}},~~\underbrace{\psi_{12},\psi_{13},\psi_{42},\psi_{43},\psi_{21},\psi_{31},\psi_{24},\psi_{34}}_{\text{indices in $E$}},~~\underbrace{\psi_{22},\psi_{23},\psi_{32} ,\psi_{33}}_{\text{indices in $D$}}~] \label{eq:order-psi-two}
\end{align}

\newpage

\begin{theorem}
\label{thm:hmt-2d}
Let $\psi_{\ell m}^I$ denote the derivative-preserving scaling of $\psi_{\ell m}$ to $I^2$, i.e.
\begin{align*}
\psi_{\ell m}^I & := \psi_{\ell}((x+1)/2)\psi_m((y+1)/2), & \ell m\in V,\\
\psi_{\ell m}^I & := 2\psi_{\ell}((x+1)/2)\psi_m((y+1)/2), & \ell m\in E.
\end{align*}
The set $\HermStwo$ has the following properties:
\begin{enumerate}
\renewcommand{\labelenumi}{(\roman{enumi})}
\item $\HermStwo$ is a basis for $\cS_3(I^2)$.
\item For any $\ell m\in V\cup E$,\quad $\xi_{\ell m}$ is identical to $\psi_{\ell m}^I$ on the edges of $I^2$.
\item $\HermStwo$ is a geometric decomposition of $\cS_3(I^2)$.
\end{enumerate}
\end{theorem}

\begin{proof}
\smartqed
The proof follows that of Theorem~\ref{thm:bern-2d} so we abbreviate proof details that are similar.
For (i), let $\HermStwoprime$ denote the derivative-preserving scaling of $\HermStwo$ to $[0,1]^2$; the scale factor is $1/2$ for functions with indices in $E$.
Given the basis orderings in (\ref{eq:order-vt-two})-(\ref{eq:order-psi-two}), we have
\begin{equation}
\label{eq:bi-hmtS-HhmtB}
\HermStwoprime = \Hmat\Hmttwo
\end{equation}
where $\Hmat$ is the $12\times 16$ matrix with the structure
\begin{equation}
\label{eq:H-struc}
\Hmat:= \left[ 
\begin{array}{c|c} ~\id~ & ~\Hmat'~ \end{array}
\right],
\end{equation}
where $\id$ is the $12\times 12$ identity matrix and $\Hmat'$ is the $12\times 4$ matrix with 
\begin{equation}
\label{eq:def-H-prime}
\Hmat' =
\left[
\begin{array}{rrrr}
 -1 & 1 & 1 & -1 \\
 1 & -1 & -1 & 1 \\
 1 & -1 & -1 & 1 \\
 -1 & 1 & 1 & -1 \\
 -1 & 0 & 1 & 0 \\
 0 & -1 & 0 & 1 \\
 1 & 0 & -1 & 0 \\
 0 & 1 & 0 & -1 \\
 -1 & 1 & 0 & 0 \\
 0 & 0 & -1 & 1 \\
 1 & -1 & 0 & 0 \\
 0 & 0 & 1 & -1 \\
\end{array}
\right ]
.
\end{equation}
Denote the entries of $\Hmat$ by $h^{\ell m}_{ij}$ (cf.\ (\ref{eq:B-entries})).
Recalling (\ref{eq:def-veps}), observe that for each $ij\in X$,
\begin{equation}
\label{eq:H-coeff-constr}
\veps_{r,i}\veps_{s,j} = \sum_{\ell m\in V\cup E} \veps_{r,\ell}\veps_{s,m} h^{\ell m}_{ij},
\end{equation}
for all $(r,s)$ pairs such that $\sldeg(x^ry^s)\leq 3$.
Similar to the Bernstein case, we substitute (\ref{eq:H-coeff-constr}) into (\ref{eq:bi-hmt-prec}), swap the order of summation and regroup, yielding
\begin{align*}
x^ry^s & = \sum_{\ell m\in V\cup E} {\veps_{r,\ell}}{\veps_{s,m}} \left(\sum_{ij\in X} h^{\ell m}_{ij}\psi_{ij} \right).
\end{align*}
The inner summation is exactly $\vt_{\ell m}^{[0,1]}$ by (\ref{eq:bi-hmtS-HhmtB}), implying that
\begin{equation}
\label{eq:bi-hmtS-prec}
x^ry^s = \sum_{\ell m\in V\cup E} {\veps_{r,\ell}}{\veps_{s,m}} \vt_{\ell m}^{[0,1]},
\end{equation}
for all $(r,s)$ pairs with $\sldeg(x^ry^s)\leq 3$, proving that $\HermStwoprime$ is a basis for $\cS_3([0,1]^2)$.
By derivative-preserving scaling, $\HermStwo$ is a basis for $\cS_3(I^2)$.

For (ii), observe that for any $ij\in D$, $\psi_{ij}\equiv 0$ on the edges of $[0,1]^2$ by virtue of the bicubic Hermite basis functions' definition.
By the structure of $\Hmat$ from (\ref{eq:H-struc}), we see that for any $\ell m\in V\cup E$,
\begin{equation}
\label{eq:vt-alt-expr}
\vt_{\ell m}^{[0,1]}=\psi_{\ell m}+\sum_{ij\in D}h^{\ell m}_{ij}\psi_{ij}.
\end{equation}
Thus, on the edges of $[0,1]^2$, $\vt_{\ell m}^{[0,1]}$ and $\psi_{\ell m}$ are identical.
After scaling back, we have $\vt_{\ell m}$ and $\psi_{\ell m}^I$ identical on the edges of $I^2$, as desired.

For (iii), the geometric decomposition is given by the indices of the basis functions, i.e.\ the function $\vt_{\ell m}$ is associated to the domain point for $\ell m\in V\cup E$.
This follows immediately from (ii), the fact that $\Hmttwo$ is a tensor product basis, and Remark~\ref{rmk:dofs}.
\qed
\end{proof}

%%  Hermite figure
\begin{figure}[h]
\begin{tabular}{ccc}
\parbox{.32\textwidth}{\includegraphics[width=.32\textwidth]{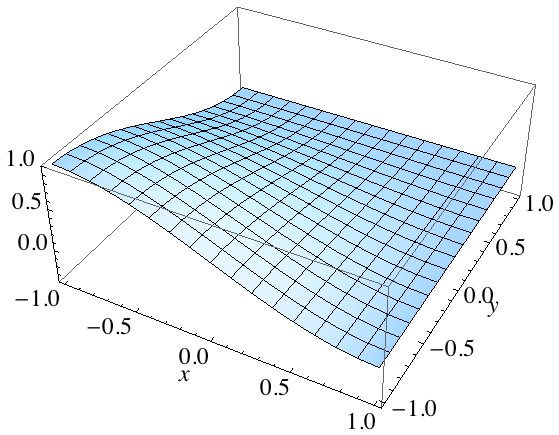}} &
\parbox{.32\textwidth}{\includegraphics[width=.32\textwidth]{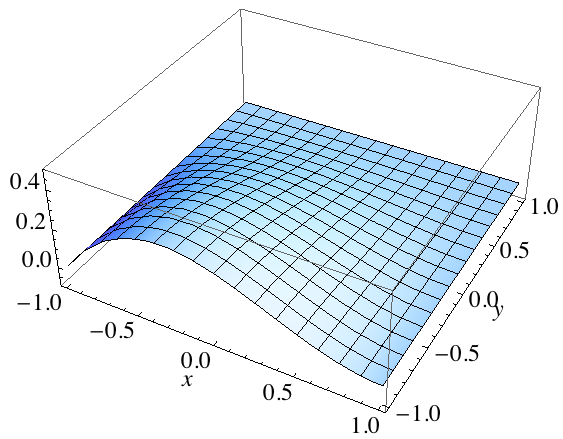}} &
\parbox{.32\textwidth}{\includegraphics[width=.32\textwidth]{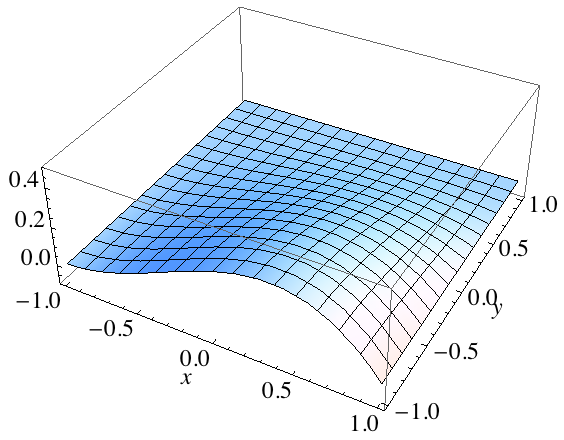}} 
\\
{$\psi_{11}^I$} & {$\psi_{21}^I$} & {$\psi_{31}^I$}
\\
\parbox{.32\textwidth}{\includegraphics[width=.32\textwidth]{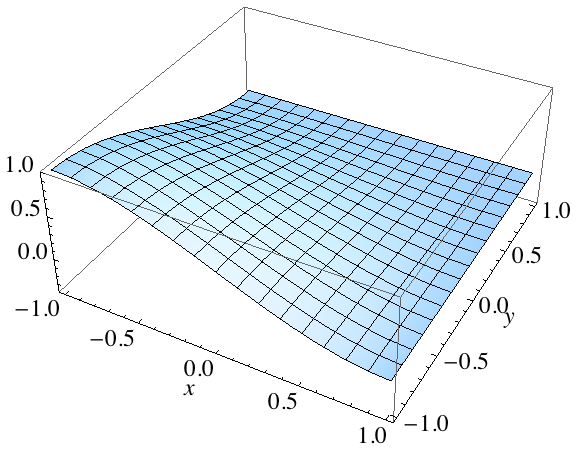}} &
\parbox{.32\textwidth}{\includegraphics[width=.32\textwidth]{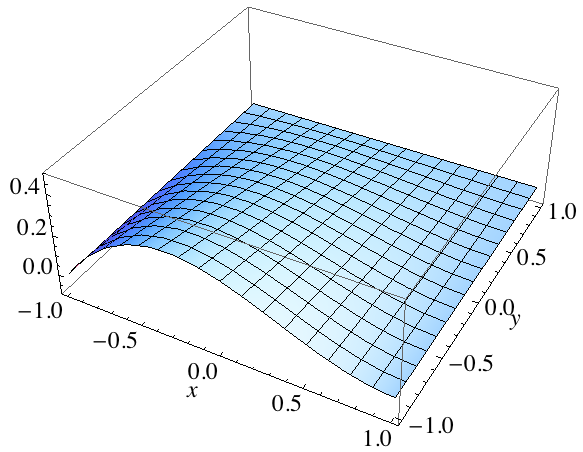}} &
\parbox{.32\textwidth}{\includegraphics[width=.32\textwidth]{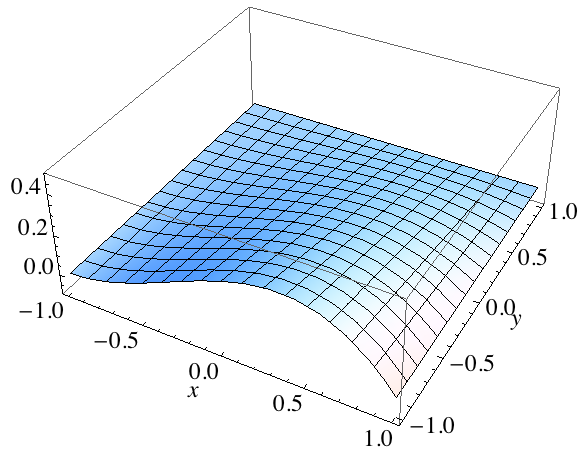}} 
\\
{$\vt_{11}$} & {$\vt_{21}$} & {$\vt_{31}$}
\end{tabular}
\caption{The top row shows 3 of the 16 bicubic Hermite functions on $I^2$ while the bottom row shows 3 of the 12 cubic Hermite style serendipity functions.  The visual differences are subtle, although some changes in concavity can be observed.  Note that functions in the same column have the same values on the edges of $I^2$.}
\label{fg:hmt-fns}
\end{figure}

\section{Local Bases for $\mathcal{S}_3(I^3)$}
\label{sec:3d-bases}

\begin{figure}[h]
\centering
\includegraphics[scale=.35]{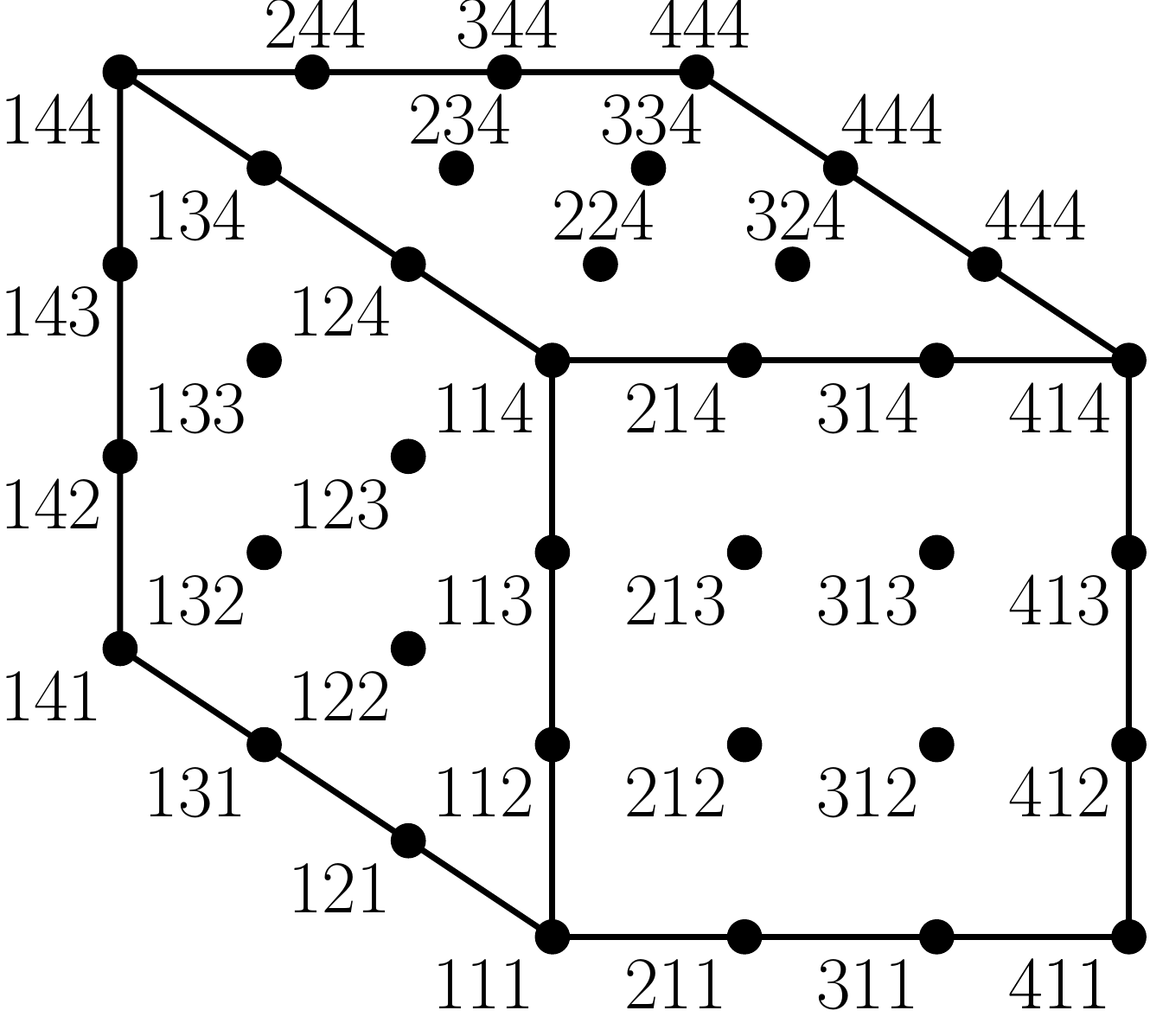}\qquad
\includegraphics[scale=.35]{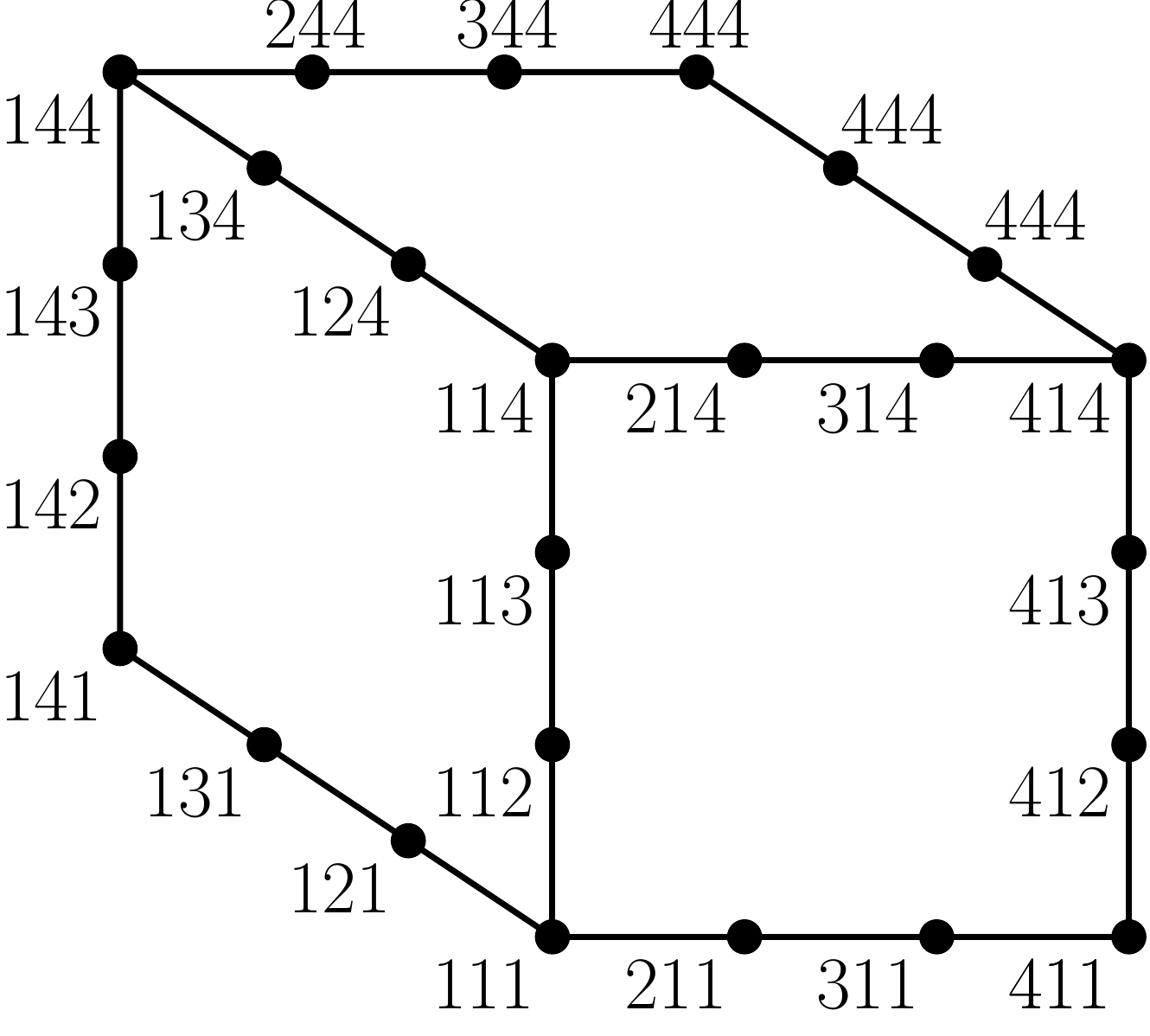}
\caption{On the left, ordered triplets from $Y$ are shown next to the domain point of $[0,1]^3$ to which they correspond.  Points hidden by the perspective are not shown.  The origin is at the point labelled 111; the positive $x$, $y$ and $z$ axes go right, back, and up, respectively.  On the right, only those indices used for the serendipity basis are shown.   The correspondences $V\leftrightarrow$ vertices, $E\leftrightarrow$ edge points, $F\leftrightarrow$ face interior points, and $M\leftrightarrow$ domain interior points are evident. }
\label{fg:cube-tricubic}
\end{figure}

Before defining local bases on the cube, we fix notation for the domain points to which they are associated.
For $[0,1]^3$, define the set of ordered triplets
\[
Y := \left\{ \{i,j,k\} ~~|~~ i,j,k\in \{1, \ldots, 4\} \right\}.
\]
Then $Y$ is the disjoint union $V\cup E\cup F\cup M$ where
\begin{align}
V & := \left\{ \{i,j,k\}\in Y ~~|~~ i,j,k\in \{1, 4\} \right\};\label{eq:cube-dom-V} \\
E & := \left\{ \{i,j,k\}\in Y ~~|~~ \text{exactly two of $i,j,k$ are elements of }\{1, 4\}\right\}; \label{eq:cube-dom-E} \\
F & := \left\{ \{i,j,k\}\in Y ~~|~~ \text{exactly one of $i,j,k$ is an element of }\{1,4\}\right\}; \label{eq:cube-dom-F} \\
M & := \left\{ \{i,j,k\}\in Y ~~|~~ i,j,k\in \{2, 3\} \right\} \label{eq:cube-dom-M} 
\end{align}
The $V$ indices are associated with vertices of $[0,1]^3$, the $E$ indices to edges of $[0,1]^3$, the $F$ indices to face interior points of $[0,1]^3$, and the $M$ vertices to the domain interior of $[0,1]^3$.
The relation between indices and domain points of the cube is shown in Figure~\ref{fg:cube-tricubic}.

\subsection{A local Bernstein style basis for $\mathcal{S}_3(I^3)$}
\label{subsec:3d-bern}

Under the notation and conventions established in Section~\ref{sec:bkgd}, we are ready to establish a local Bernstein style basis for $\mathcal{S}_3(I^3)$ where $I:=[-1,1]$.
In Figure~\ref{fig:bern-ser-3d-def}, we define a set of 32 functions, indexed by ${V\cup E\subset Y}$; note the scaling by $1/32$.
We fix the following basis orderings, with omitted basis functions ordered lexicographically by index.
\begin{align}
\BernSthree &:= [~\underbrace{\xi_{111},\ldots,\xi_{444}}_{\text{indices in $V$}},~~\underbrace{\xi_{112},\ldots,\xi_{443}}_{\text{indices in $E$}}~],\label{eq:order-xi-three}\\
\Bern & := [~
\underbrace{\beta_{111},\ldots,\beta_{444}}_{\text{indices in $V$}},~~
\underbrace{\beta_{112},\ldots,\beta_{443}}_{\text{indices in $E$}},~~
\underbrace{\beta_{122},\ldots,\beta_{433}}_{\text{indices in $F$}},~~
\underbrace{\beta_{222},\ldots,\beta_{333}}_{\text{indices in $M$}},
~] 
\label{eq:order-beta-three}
\end{align}

\begin{figure}
\begin{align*}
\BernSthree =
\begin{bmatrix}
\xi_{111} \\
\xi_{114} \\
\xi_{141}  \\
\xi_{144} \\
\xi_{411} \\
\xi_{414} \\
\xi_{441} \\
\xi_{444} \\
\xi_{112} \\
\xi_{113} \\
\xi_{121} \\
\xi_{124} \\
\xi_{131} \\
\xi_{134} \\
\xi_{142} \\
\xi_{143} \\
\xi_{211} \\
\xi_{214} \\
\xi_{241} \\
\xi_{244} \\
\xi_{311} \\
\xi_{314} \\
\xi_{341} \\
\xi_{344} \\
\xi_{412} \\
\xi_{413} \\
\xi_{421} \\
\xi_{424} \\
\xi_{431} \\
\xi_{434} \\
\xi_{442} \\
\xi_{443}
\end{bmatrix}
 & =
\begin{bmatrix}
(1-x) (1-y) (1-z) (-5-2 x+x^2-2 y+y^2-2 z+z^2) \\
(1-x) (1-y) (z+1) (-5-2 x+x^2-2 y+y^2+2 z+z^2) \\
(1-x) (y+1) (1-z) (-5-2 x+x^2+2 y+y^2-2 z+z^2) \\
(1-x) (y+1) (z+1) (-5-2 x+x^2+2 y+y^2+2 z+z^2) \\
(x+1) (1-y) (1-z) (-5+2 x+x^2-2 y+y^2-2 z+z^2) \\
(x+1) (1-y) (z+1) (-5+2 x+x^2-2 y+y^2+2 z+z^2) \\
(x+1) (y+1) (1-z) (-5+2 x+x^2+2 y+y^2-2 z+z^2) \\
(x+1) (y+1) (z+1) (-5+2 x+x^2+2 y+y^2+2 z+z^2) \\
(1-x) (1-y) (1-z)^2 (z+1) \\
(1-x) (1-y) (1-z) (z+1)^2 \\
(1-x) (1-y)^2 (y+1) (1-z) \\
(1-x) (1-y)^2 (y+1) (z+1) \\
(1-x) (1-y) (y+1)^2 (1-z) \\
(1-x) (1-y) (y+1)^2 (z+1) \\
(1-x) (y+1) (1-z)^2 (z+1) \\
(1-x) (y+1) (1-z) (z+1)^2 \\
(1-x)^2 (x+1) (1-y) (1-z) \\
(1-x)^2 (x+1) (1-y) (z+1) \\
(1-x)^2 (x+1) (y+1) (1-z) \\
(1-x)^2 (x+1) (y+1) (z+1) \\
(1-x) (x+1)^2 (1-y) (1-z) \\
(1-x) (x+1)^2 (1-y) (z+1) \\
(1-x) (x+1)^2 (y+1) (1-z) \\
(1-x) (x+1)^2 (y+1) (z+1) \\
(x+1) (1-y) (1-z)^2 (z+1) \\
(x+1) (1-y) (1-z) (z+1)^2 \\
(x+1) (1-y)^2 (y+1) (1-z) \\
(x+1) (1-y)^2 (y+1) (z+1) \\
(x+1) (1-y) (y+1)^2 (1-z) \\
(x+1) (1-y) (y+1)^2 (z+1) \\
(x+1) (y+1) (1-z)^2 (z+1) \\
(x+1) (y+1) (1-z) (z+1)^2  
\end{bmatrix}
\cdot \frac 1{32}
\end{align*}
\caption{Bernstein-style basis functions for $\cS_3(I^3)$ with properties given by Theorem~\ref{thm:bern-3d}}
\label{fig:bern-ser-3d-def}
\end{figure}

\begin{theorem}
\label{thm:bern-3d}
Let $\beta_{\ell mn}^I$ denote the scaling of $\beta_{\ell mn}$ to $I^3$, i.e.\
\[\beta_{\ell mn}^I :=\beta_{\ell}((x+1)/2)\beta_m((y+1)/2)\beta_n((z+1)/2).\]
The set $\BernSthree$ has the following properties:
\begin{enumerate}
\renewcommand{\labelenumi}{(\roman{enumi})}
\item $\BernSthree$ is a basis for $\cS_3(I^3)$.
\item $\BernSthree$ reduces to $\BernStwo$ on faces of $I^3$.
\item For any $\ell mn\in V\cup E$,\quad $\xi_{\ell mn}$ is identical to $\beta_{\ell mn}^I$ on edges of $I^3$.
\item $\BernSthree$ is a geometric decomposition of $\cS_3(I^3)$.
\end{enumerate}
\end{theorem}

\begin{proof}
\smartqed
The proof is similar to that of Theorem~\ref{thm:bern-2d}.
Note that for (ii), the claim can be confirmed directly by calculation, for instance, $\xi_{111}(x,y,-1)=\xi_{11}(x,y)$ or $\xi_{142}(x,1,z)=\xi_{12}(x,z)$.
The proof is similar to that of Theorem~\ref{thm:bern-2d} so we abbreviate some details.
For (i), let $\BernSthreeprime$ denote the scaling of $\BernSthree$ to $[0,1]^3$.
Given the basis orderings in (\ref{eq:order-xi-three})-(\ref{eq:order-beta-three}), we have
\begin{equation}
\label{eq:tri-bernS-AbernB}
\BernSthreeprime = \U\BernThree
\end{equation}
where $\U$ is the $32\times 64$ matrix with the structure
\begin{equation}
\label{eq:K-struc}
\U:= \left[ 
\begin{array}{c|c} ~\id ~& ~\U'~ \end{array}
\right],
\end{equation}
where $\id$ is the $32\times 32$ identity matrix and $\U'$ is a specific $32\times 32$ matrix whose entries are integers ranging from -16 to 4.
Instead of writing out $\U'$ in its entirety, we describe its properties and how it can be constructed (cf.\ Remark~\ref{rmk:proof-tech}).

Using $ijk\in Y$ to denote an index for $\beta_{ijk}$ and $\ell mn\in V\cup E\subset Y$ to denote an index for $\xi_{\ell mn}^{[0,1]}$, the entries of $\U$ will be denoted by $k^{\ell mn}_{ijk}$ so that
\begin{equation}
\label{eq:K-entries}
\U:= \left[ 
\begin{array}{ccccc}
u^{111}_{111} & \cdots & u^{111}_{ijk} & \cdots &  u^{111}_{333} \\ 
\vdots & \ddots & \vdots & \ddots  & \vdots \\
u^{\ell mn}_{111} & \cdots & u^{\ell mn}_{ijk} & \cdots &  u^{\ell mn}_{333} \\ 
\vdots & \ddots & \vdots & \ddots  & \vdots \\
u^{443}_{111} & \cdots & u^{443}_{ijk} & \cdots &  u^{443}_{333}
\end{array}
\right].
\end{equation}
The columns of $\U$ satisfy the relationship
\begin{equation}
\label{eq:K-coeff-constr}
{3-r\choose 4-i}{3-s\choose 4-j}{3-t\choose 4-k} = \sum_{\ell mn\in V\cup E} {3-r\choose 4-\ell}{3-s\choose 4-m}{3-t\choose 4-n} u^{\ell mn}_{ijk},
\end{equation}
for all $(r,s,t)$ tuples such that  $\sldeg(x^ry^sz^t)\leq 3$.
This property defines the entries of $\U$ uniquely since for each $ijk\in Y$ it gives an invertible linear system of 32 equations with 32 unknowns whose solution is the $ijk$ column of $\U$.
The $(r,s,t)$ tuples should be taken in the order given in (\ref{def:p3-3d})-(\ref{def:s3-3d}).
See Remark~\ref{rmk:proof-tech} and the text after (\ref{eq:A-coeff-constr}) for more on this process.

As in previous proofs, regrouping and recognizing an expression for $\xi_{\ell mn}^{[0,1]}$ gives
\begin{equation}
\label{eq:tri-bernS-prec}
x^ry^sz^t = \sum_{\ell mn\in V\cup E} {3-r\choose 4-\ell}{3-s\choose 4-m}{3-t\choose 4-n}\xi_{\ell mn},
\end{equation}
proving, after scaling, that $\BernSthree$ is a basis for $\cS_3(I^3)$.

For (ii), the claim can be confirmed directly by calculation, e.g.\ $\xi_{111}(x,y,-1)=\xi_{11}(x,y)$ or $\xi_{142}(x,1,z)=\xi_{12}(x,z)$, etc.

For (iii), note that an edge of $[0,1]^3$ is described by two equations of the form $\{x,y,~\text{or}~z\}=\{0~\text{or}~1\}$ where two distinct variables must be chosen for the two equations.
Since $\beta_2(t)$ and $\beta_3(t)$ are equal to 0 at $t=0$ and $t=1$, $\beta_{ijk}\equiv 0$ on the edges of $[0,1]^2$ for any $ijk\in M$.
Further, for $ijk\in F$, without loss of generality, assume that $i\in\{1,4\}$ so that $j,k\in\{2,3\}$.
Since every edge is described by at least one equation of the form $\{y~\text{or}~z\}=\{0~\text{or}~1\}$, either $\beta_j(y)$ or $\beta_k(z)$ is identically zero on every edge.
Thus, for $ijk\in F\cup M$, $\beta_{ijk}\equiv 0$ on the edges of $[0,1]^3$.

By the structure of $\U$ from (\ref{eq:K-struc}), we see that for any $\ell mn\in V\cup E$,
\begin{equation}
\label{eq:xi-alt-expr-3d}
\xi_{\ell mn}^{[0,1]}=\beta_{\ell mn}+\sum_{ijk\in F\cup M}u^{\ell mn}_{ijk}\beta_{ijk}.
\end{equation}
Thus, on the edges of $[0,1]^3$, $\xi_{\ell mn}^{[0,1]}$ and $\beta_{\ell mn}$ are identical.
After scaling back, we have $\xi_{\ell mn}$ and $\beta_{\ell mn}^I$ identical on the edges of $I^3$, as desired.

For (iv), the geometric decomposition is given by the indices of the basis functions, i.e.\ the function $\xi_{\ell mn}$ is associated to the domain point for $\ell mn\in V\cup E$.
This follows immediately from (ii) and (iii), the fact that $\BernThree$ is a tensor product basis, and Remark~\ref{rmk:dofs}.
\qed
\end{proof}

\subsection{A local Hermite style basis for $\mathcal{S}_3(I^3)$}
\label{subsec:3d-hmt}

We now establish a local Hermite style basis $\HermSthree$ for $\cS_3(I^3)$ using the tricubic Hermite basis $\Hmtthree$ for $\cQ_3([0,1]^3)$.
In Figure~\ref{fig:hmt-ser-3d-def}, we define a set of 32 functions, indexed by ${V\cup E\subset Y}$; note the scaling by $1/16$.
We fix the following basis orderings, with omitted basis functions ordered lexicographically by index.
\begin{align}
\HermSthree &:= [~\underbrace{\vt_{111},\ldots,\vt_{444}}_{\text{indices in $V$}},~~\underbrace{\vt_{112},\ldots,\vt_{443}}_{\text{indices in $E$}}~],\label{eq:order-vt-three}\\
\Bern & := [~
\underbrace{\psi_{111},\ldots,\psi_{444}}_{\text{indices in $V$}},~~
\underbrace{\psi_{112},\ldots,\psi_{443}}_{\text{indices in $E$}},~~
\underbrace{\psi_{122},\ldots,\psi_{433}}_{\text{indices in $F$}},~~
\underbrace{\psi_{222},\ldots,\psi_{333}}_{\text{indices in $M$}},
~] \label{eq:order-psi-three}
\end{align}

\begin{figure}
\begin{align*}
\HermSthree =
\begin{bmatrix}
\vt_{111} \\
\vt_{114} \\
\vt_{141}  \\
\vt_{144} \\
\vt_{411} \\
\vt_{414} \\
\vt_{441} \\
\vt_{444} \\
\vt_{112} \\
\vt_{113} \\
\vt_{121} \\
\vt_{124} \\
\vt_{131} \\
\vt_{134} \\
\vt_{142} \\
\vt_{143} \\
\vt_{211} \\
\vt_{214} \\
\vt_{241} \\
\vt_{244} \\
\vt_{311} \\
\vt_{314} \\
\vt_{341} \\
\vt_{344} \\
\vt_{412} \\
\vt_{413} \\
\vt_{421} \\
\vt_{424} \\
\vt_{431} \\
\vt_{434} \\
\vt_{442} \\
\vt_{443}
\end{bmatrix}
 & =
\begin{bmatrix}
-(1-x) (1-y) (1-z) (-2+x+x^2+y+y^2+z+z^2) \\
-(1-x) (1-y) (z+1) (-2+x+x^2+y+y^2-z+z^2) \\
-(1-x) (y+1) (1-z) (-2+x+x^2-y+y^2+z+z^2) \\
-(1-x) (y+1) (z+1) (-2+x+x^2-y+y^2-z+z^2) \\
-(x+1) (1-y) (1-z) (-2-x+x^2+y+y^2+z+z^2) \\
-(x+1) (1-y) (z+1) (-2-x+x^2+y+y^2-z+z^2) \\
-(x+1) (y+1) (1-z) (-2-x+x^2-y+y^2+z+z^2) \\
-(x+1) (y+1) (z+1) (-2-x+x^2-y+y^2-z+z^2) \\
 (1-x) (1-y) (1-z)^2 (z+1) \\
 (1-x) (1-y) (1-z) (z+1)^2 \\
 (1-x) (1-y)^2 (y+1) (1-z) \\
 (1-x) (1-y)^2 (y+1) (z+1) \\
 (1-x) (1-y) (y+1)^2 (1-z) \\
 (1-x) (1-y) (y+1)^2 (z+1) \\
 (1-x) (y+1) (1-z)^2 (z+1) \\
 (1-x) (y+1) (1-z) (z+1)^2 \\
 (1-x)^2 (x+1) (1-y) (1-z) \\
 (1-x)^2 (x+1) (1-y) (z+1) \\
 (1-x)^2 (x+1) (y+1) (1-z) \\
 (1-x)^2 (x+1) (y+1) (z+1) \\
 (1-x) (x+1)^2 (1-y) (1-z) \\
 (1-x) (x+1)^2 (1-y) (z+1) \\
 (1-x) (x+1)^2 (y+1) (1-z) \\
 (1-x) (x+1)^2 (y+1) (z+1) \\
 (x+1) (1-y) (1-z)^2 (z+1) \\
 (x+1) (1-y) (1-z) (z+1)^2 \\
 (x+1) (1-y)^2 (y+1) (1-z) \\
 (x+1) (1-y)^2 (y+1) (z+1) \\
 (x+1) (1-y) (y+1)^2 (1-z) \\
 (x+1) (1-y) (y+1)^2 (z+1) \\
 (x+1) (y+1) (1-z)^2 (z+1) \\
 (x+1) (y+1) (1-z) (z+1)^2 
\end{bmatrix}
\cdot \frac 1{16}
\label{eq:hmt-ser-3d-def}
\end{align*}
\caption{Hermite-style basis functions for $\cS_3(I^3)$ with properties given by Theorem~\ref{thm:hmt-3d}}
\label{fig:hmt-ser-3d-def}
\end{figure}

\begin{theorem}
\label{thm:hmt-3d}
Let $\psi_{\ell mn}^I$ denote the derivative-preserving scaling of $\psi_{\ell mn}$ to $I^3$, i.e.
\begin{align*}
\psi_{\ell m}^I & := \psi_{\ell}((x+1)/2)\psi_m((y+1)/2)\psi_n((z+1)/2), & \ell mn\in V,\\
\psi_{\ell mn}^I & := 2\psi_{\ell}((x+1)/2)\psi_m((y+1)/2)\psi_n((z+1)/2), & \ell mn\in E.
\end{align*}
The set $\HermSthree$ has the following properties:
\begin{enumerate}
\renewcommand{\labelenumi}{(\roman{enumi})}
\item $\HermSthree$ is a basis for $\cS_3(I^3)$.
\item $\HermSthree$ reduces to $\HermStwo$ on faces of $I^3$.
\item For any $\ell mn\in V\cup E$,\quad $\vt_{\ell mn}$ is identical to $\psi_{\ell mn}^I$ on edges of $I^3$.
\item $\HermSthree$ is a geometric decomposition of $\cS_3(I^3)$.
\end{enumerate}
\end{theorem}

The proof is similar to that of Theorem~\ref{thm:hmt-2d}.

\begin{proof}
\smartqed
The proof is similar to that of Theorem~\ref{thm:hmt-2d} so we abbreviate some details.
For (i), let $\HermSthreeprime$ denote the scaling of $\HermSthree$ to $[0,1]^3$.
Given the basis orderings in (\ref{eq:order-vt-three})-(\ref{eq:order-psi-three}), we have
\begin{equation}
\label{eq:tri-hmtS-AhmtB}
\HermSthreeprime = \W\Hmtthree
\end{equation}
where $\W$ is the $32\times 64$ matrix with the structure
\begin{equation}
\label{eq:W-struc}
\W := \left[ 
\begin{array}{c|c} ~\id~ & ~\W'~ \end{array}
\right],
\end{equation}
where $\id$ is the $32\times 32$ identity matrix and $\W'$ is a specific $32\times 32$ matrix whose entries are in $\{-1,0,1\}$.
The matrix $\W$ is constructed similarly to the matrix $\U$ from the proof of Theorem~\ref{thm:bern-3d}; the columns of $\W$ satisfy the relationship
\begin{equation}
\label{eq:W-coeff-constr}
\veps_{r,i}\veps_{s,j}\veps_{t,k} = \sum_{\ell mn\in V\cup E} \veps_{r,\ell}\veps_{s,m}\veps_{t,n}w^{\ell mn}_{ijk},
\end{equation}
for all $(r,s,t)$ tuples such that  $\sldeg(x^ry^sz^t)\leq 3$.
Similar to previous proofs, this yields
\begin{equation}
\label{eq:tri-hmtS-prec}
x^ry^sz^t = \sum_{\ell mn\in V\cup E} \veps_{r,\ell}\veps_{s,m}\veps_{t,n} \vt_{\ell mn},
\end{equation}
proving, after derivative-preserving scaling, that $\HermSthree$ is a basis for $\cS_3(I^3)$.

For (ii)-(iv), the proof is similar to the corresponding parts of the proof of Theorem~\ref{thm:bern-3d}.
\qed
\end{proof}

\section{Conclusions and Future Directions}
\label{sec:conc}

The basis functions presented in this work are well-suited for use in finite element applications, as discussed in the introduction.
For geometric modeling purposes, some adaptation of traditional techniques will be required as the bases do not have the classical properties of positivity and do not form a partition of unity.
Nevertheless, we are already witnessing the successful implementation of the basis $\HermSthree$ in the geometric modeling and finite element analysis package Continuity developed by the Cardiac Mechanics Research Group at UC San Diego.
In that context, the close similarities of $\HermSthree$ and $[\psi^3]$ has allowed a straightforward implementation procedure with only minor adjustments to the geometric modeling subroutines.

Additionally, the proof techniques used for the theorems suggest a number of promising extensions.
Similar techniques should be able to produce Bernstein style bases for higher polynomial order serendipity spaces, although the introduction of interior degrees of freedom that occurs when $r>3$ requires some additional care to resolve.
Some higher order Hermite style bases may also be available, although the association of directional derivative values to vertices is somewhat unique to the $r=3$ case.
Pre-conditioners for finite element methods employing our bases are still needed, as is a thorough analysis of the tradeoffs between the approach outlined here and alternative approaches to basis reduction, such as static condensation.
The fact that all the functions defined here are fixed linear combinations of standard bicubic or tricubic basis functions suggests that appropriate pre-conditioners will have a straightforward and computationally advantageous construction.

\begin{acknowledgement}
Support for this work was provided in part by NSF Award~0715146 and the National Biomedical Computation Resource while the author was at the University of California, San Diego.
\end{acknowledgement}
%
%\section*{Appendix}
%\addcontentsline{toc}{section}{Appendix}
%
%
\bibliographystyle{abbrv}
\bibliography{akg}

\end{document}